\documentclass [12pt]{article}
\title {A note on the $NFI$-topology}
\author {Ziv Shami\\Ariel University}

\newtheorem {theorem}{Theorem}[section]

\newtheorem {definition}[theorem]{Definition}

\newtheorem {fact}[theorem]{Fact}

\newtheorem {remark}[theorem]{Remark}

\newtheorem {proposition}[theorem]{Proposition}

\newtheorem {claim}[theorem]{Claim}

\def\proof {\noindent \textbf{Proof:} }

\def\qed {$\ \ \ \ \Box$}

\newsavebox{\indbin}
\savebox{\indbin}{\begin{picture}(0,0)
\newlength{\gnu}
\settowidth{\gnu}{$\smile$} \setlength{\unitlength}{.5\gnu} \put(-1,-.65){$\smile$}
\put(-.25,.1){$|$}
\end{picture}}
\newcommand{\nonfork}[3]
{\mbox{$\begin{array}{ccc} \mbox{$#1$} & \usebox{\indbin} & \mbox{$#2$} \\
        & \mbox{$#3$} &
\end{array}$}}
\newcommand{\nonforkempty}[2]
{\mbox{$\begin{array}{ccc} \mbox{$#1$} & \usebox{\indbin} & \mbox{$#2$}
\end{array}$}}

\newsavebox{\sindbin}
\savebox{\sindbin}{\begin{picture}(0,0)
\newlength{\sgnu}
\settowidth{\sgnu}{$\smile$} \setlength{\unitlength}{.5\sgnu} \put(-1,-.65){$\smile$}
\put(-.25,.1){$|s$}
\end{picture}}

\newsavebox{\starindbin}
\savebox{\starindbin}{\begin{picture}(0,0)
\newlength{\stargnu}
\settowidth{\stargnu}{$\smile$} \setlength{\unitlength}{.5\stargnu} \put(-1,-.65){$\smile$}
\put(-.25,.1){$|*$}
\end{picture}}

\newsavebox{\qindbin}
\savebox{\qindbin}{\begin{picture}(0,0)
\newlength{\qgnu}
\settowidth{\qgnu}{$\smile$} \setlength{\unitlength}{.5\qgnu} \put(-1,-.65){$\smile$}
\put(-.25,.1){$|_{qf}$}
\end{picture}}

\newsavebox{\minusindbin}
\savebox{\minusindbin}{\begin{picture}(0,0)
\newlength{\minusgnu}
\settowidth{\minusgnu}{$\smile$} \setlength{\unitlength}{.5\minusgnu} \put(-1,-.65){$\smile$}
\put(-.25,.1){$|{^-}$}
\end{picture}}
\newcommand{\minusnonfork}[3]
{\mbox{$\begin{array}{ccc} \mbox{$#1$} & \usebox{\minusindbin} & \mbox{$#2$} \\
        & \mbox{$#3$} &
\end{array}$}}
\newcommand{\minusnonforkempty}[2]
{\mbox{$\begin{array}{ccc} \mbox{$#1$} & \usebox{\minusindbin} & \mbox{$#2$}
\end{array}$}}

\def\card #1 {{\vert #1 \vert}}

\def\CC {{\cal C}}

\def\UU {{\cal U}}

\def\FF {{\cal F}}
\def\OO {{\cal O}}
\def\MM {{\cal M}}
\def\Tp {\vert T\vert^+}

\usepackage{amsfonts}

\begin{document}
\maketitle

\begin{abstract}
The $NFI$-topology, introduced in [S0], is a topology on the Stone space of a theory $T$ that depends on a reduct $T^-$ of $T$. This topology has been used in [S0]
to describe the set of universal transducers for $(T,T^-)$ (invariants sets that translates forking-open sets in $T^-$ to forking-open sets in $T$).
In this paper we show that in contrast to the stable case, the $NFI$-topology need not be invariant over parameters in $T^-$ but a weak version of this holds for any simple $T$.
We also note that for the lovely pair expansions, of theories with the \em wnfcp\em , the topology is invariant over $\emptyset$ in $T^-$.

\end{abstract}


\section{Introduction}

Recall the definition of the forking topology for a theory $T$ [S1]: if  $A\subseteq\CC$ is small and  $x$ is a finite tuple of variable,  a set $U\subseteq S_x(A)$ is said to be \em a basic forking-open set over $A$ for $T$ \em if there exists $\phi(x,y)\in L(A)$ such that $$U=U^T_\phi\equiv\{p\in S_x(A) \vert\ \phi(a,y)\ L-\mbox{forks\ over}\ A\ \mbox{for\ all\ } a\models p \}.$$
The family of basic forking-open sets over $A$ is closed under finite intersections and thus form a basis for a unique topology on $S_x(A)$. Note that the forking-topology on $S_x(A)$ refines the Stone-topology.
For a simple $L$-theory $T$ and a reduct $T^-$ to a sublanguage $L^-\subseteq L$,  an $\emptyset$-invariant set $\Gamma(x)$ in a monster $\CC$ of $T$ is called a \em universal transducer \em if  for every formula $\phi^-(x,y)\in L^-$ and every $a$,  $$\phi^-(x,a)\ L^-\mbox{-forks\ over}\ \emptyset\ \mbox{\ iff\ } \Gamma(x)\wedge \phi^-(x,a)\ L\mbox{-forks\ over}\ \emptyset.$$ Moreover, there is a greatest universal transducer $\tilde \Gamma_x$ and it is type-definable. In particular, the forking topology for $T$ on $S_y(T)$ refines the forking topology for $T^-$ on $S_y(T^-)$ for all $y$ [S0]. In  [S0] a new topology on $S_x(T)$ is introduced  (the $NFI$-topology, see section 2) such that an $\emptyset$-invariant set $\Gamma(x)$ in $\CC$ is a universal transducer iff it is a dense subset of $\tilde\Gamma_x$ in the relative $NFI$-topology on $\tilde\Gamma_x$. If $T$ is stable, a subset of $S_x(T)$ is open in the $NFI$-topology iff it is a union of $L$-definable sets over $\emptyset$ that are $L^-$-definable with parameters [S0, Lemma 2.20].
The main goal of this paper is examine to what extent the latest property can be generalized to the simple case.  First we show that forking-open sets for $T$ defined by formulas in $L^-$ are forking-open for $T^-$
over parameters; this is, a very weak version of $L^-$-invariance of the $NFI$-topology of $T$ over parameters.  Then, we show an example of a simple theory and a basic $NFI$-open set that is not invariant over parameters in the reduct. Finally, we show that for the lovely-pair expansions of theories with the \em wnfcp\em, every basic $NFI$-open set  in the expansion is in fact a basic  $NFI$-open set in the reduct and thus in particular, the $NFI$-topology of $T$ is $L^-$-invariant over parameters.

We assume basic knowledge of simple theories as in [K],[KP],[HKP]. A good textbook on simple theories is [W]. In this paper, unless otherwise stated,  $T$ will denote a complete first-order simple theory in an arbitrary language $L$ (unless otherwise stated) and we work in a $\lambda$-big model $\CC$ of $T$ (i.e. a model with the property that any expansion of it by less than $\lambda$ constants is splendid) for some large $\lambda$. We call $\CC$ the monster model. Note that any $\lambda$-big model (of any theory) is $\lambda$-saturated and $\lambda$-strongly homogeneous and that $\lambda$-bigness is preserved under reducts (by Robinson consistency theorem).

\section{Forking invariance in a reduct}
In this section $T$ denotes a simple $L$-theory and $T^-$ denotes a reduct of $T$ to a sublanguage $L^-$ of $L$ and $\CC^-=\CC\vert L^-$. As mentioned in the introduction, we know that both $\CC$ and $\CC^-$ are highly saturated and highly strongly-homogeneous. We use $\nonforkempty {}{}$ to denote independence in $\CC$, and $\minusnonforkempty{}{}$ to denote independence in $\CC^-$.
For a small set $A\subseteq \CC^{heq}$, $BDD(A)$ denotes the set of countable (length)  hyperimaginaries in $\CC^{heq}$  that are in the bounded closure of $A$ in the sense of $\CC$.
$LSTP(a)$ denotes the Lascar of $a$ in $\CC$.\\

\noindent Recall the definition of $NFI$-topology (restricted version):

\begin{definition}\label{def. NFI-top}\em
Given a finite tuple of variables $y$, a set $U=U(y)$ is \em a basic open set in the $NFI$-topology \em  on $S_y(T)$ iff
there exists a type $p(x)\in S_x(T)$  and $\phi^-(x,y)\in L^-$ such that $$U=U_{p,\phi^-}=\{b\vert\ p(x)\wedge\phi^-(x,b)\ L\mbox{-doesn't fork\ over}\ \emptyset\}.$$
\end{definition}

We introduce now a related notion:

\begin{definition}\em
1)  Let $p\in S_x(T)$ and let $\phi^-(x,y)\in L^-$. We say that \em $p$ is $L^--IF$ for $\phi^-(x,y)$ \em if the set $U_{p,\phi^-}$ is an $L^-$-invariant set over some small set.\\
2) Let $A$ be a small set. We say that a type $p\in S_x(T)$ is \em $L^--IF$ over $A$ \em if for every $\phi^-(x,y)\in L^-$, the set $U_{p,\phi^-}$ is an $L^-$-invariant set over $A$.\\
3) We say that a type $p\in S_x(T)$ is \em $L^--IF$ \em if it is  $L^--IF$ over some set $A$.
\end{definition}

\begin{remark}
Let $T$ be any complete theory. Let $\phi^-(x,y)\in L^-$. If every $p\in S_x(T)$ is $L^--IF$ for $\phi^-(x,y)$ then $U^T_{\phi^-}$ is an $L^-$-invariant set over parameters.
\end{remark}

\begin{proposition}
Let $\phi^-(x,y)\in L^-$ and let $U^T_{\phi^-}$ be the corresponding basic forking open set over $\emptyset$ for $T$.  Then $U^T_{\phi^-}$ is forking open set over some small set for $T^-$.
\end{proposition}

\proof  Let $NF_{\phi^-}^L=\{b \vert \ \phi^-(x,b)\ L-doesn't\ fork\ over\ \emptyset\}$ be the complement of  $U^T_{\phi^-}$. Working in $\CC$, let $\FF=\{J_i \vert i<\beta\}$, where $J_i=(J_{i,j} \vert  j<\Tp)$,  be a family of independent Morley sequences over $\emptyset$ in the sort of $x$, of length $\vert T\vert^+$ such that for every possible Lascar strong type of such a sequence (over $\emptyset$) there are sufficiently many independent realizations of that Lascar strong type. Let $I_{\FF}=(\bigwedge_i {J_{i,j}} \vert j<\Tp)$. It will be sufficient to prove the following.

\begin{claim}\label{claim1}
For every $b$, we have $b\in NF_{\phi^-}^L$ iff $\phi^-(x,b)$ $L^-$-doesn't fork over ${I_{\FF}}^{>\alpha}$ for all $\alpha<\vert T\vert^+$.
\end{claim}

\proof Let $b$ be such that $\phi^-(x,b)$ $L^-$-doesn't fork over ${I_{\FF}}^{>\alpha}$ for all $\alpha<\vert T\vert^+$. Then, for some $\alpha^*<\vert T\vert^+$
we have $\nonforkempty{b}{{I_{\FF}}^{>\alpha^*}}$  (otherwise, we can construct by induction a sequence $(\bar e_i\vert i<\vert T\vert^+)$ of pairwise disjoint finite subsequences of ${I_{\FF}}$ such that
$e_i$ and $b$ are $L$-dependent for every $i$ which contradicts simplicity of $T$). By the assumption, there exists $a$ such that $a\models\phi^-(a,b)$ and $\minusnonfork{a}{b}{{I_{\FF}}^{>\alpha^*}}$.
Let $a'$ realize $tp_{L^-}(a/b{I_{\FF}}^{>\alpha^*})$ such that $\nonfork{a'}{b}{{I_{\FF}}^{>\alpha^*}}$  . By the choice of $\alpha^*$ , we conclude that $\phi^-(a',b)$ and $\nonforkempty{a'}{b}$.

For the other direction, assume $\phi^-(x,b)$ $L$-doesn't fork over $\emptyset$. Let $a\models \phi^-(x,b)$ such that $\nonforkempty{a}{b}$ and let $I=(a_i\vert i<\vert T\vert^+)$ be an $L$-Morley sequence of $tp_L(a/b)$. Let $I_*\in\FF$ be such that $LSTP(I_*)=LSTP(I)$ and  $\nonforkempty{I_*}{I}$. It will be sufficient to show that $\phi^-(x,b)$ $L^-$-doesn't fork over $I_*^{>\alpha}$ for all $\alpha<\vert T\vert^+$. Assume by contradiction that $\phi^-(x,b)$ $L^-$-forks over $I_*^{>\alpha_0}$ for some $\alpha_0<\vert T\vert^+$. Let $\OO=\Tp\frown \Tp_{rev}$ be the concatenation of the $\in$-order of $\Tp$ with its reverse order. We use the notation $\Tp_{rev}=\{ m_{-i} \ \vert\ i<\Tp\}$ (where $m_{-i}$ denotes the $i$-th element from the top of $\OO$).
Let $\tilde I=(a_s \vert s\in\OO)$ be an indiscernible sequence that extends $I$. Let $\bar x=(x_i \vert i<\Tp)$ and let $\bar x_{rev}= (x_s\ \vert\ s\in \Tp_{rev})$.
Now, let $p(\bar x_{rev},I)=tp_L((a_{s} \vert s\in\Tp_{rev})/I)$ and let $p(\bar x_{rev},I_*)$ be the corresponding $BDD(\emptyset)$-conjugate of $p(\bar x_{rev},I)$. By the independence theorem in $\CC$ the type $p(\bar x_{rev},I_*)\wedge p(\bar x_{rev},I)$ $L$-doesn't fork over $\emptyset$ and in particular there exists a sequence $I'=(a'_s \vert s\in\Tp_{rev} )$ such that both $I_*\frown I'$ and $I\frown I'$ are $L$-indiscernible over $\emptyset$. Let $I'_\alpha=(a'_{m_{-(\omega\alpha+i)}}\vert i<\omega)$. Note that the sequence $(I'_\alpha \vert \alpha<\Tp)$ is an $L^-$-Morley sequence over $I_*^{>\alpha_0}$ (as ${tp_{L^-}}(I'_\alpha/(\bigcup_{\beta<\alpha} I'_\beta)\cup I_*^{>\alpha_0})$ is finitely satisfiable in $I_*^{>\alpha_0}$), hence there exists $\alpha^*<\Tp$ such that $\minusnonfork{I'_{\alpha^*}}{b}{I_*^{>\alpha_0}}$. Therefore, $\phi^-(x,b)$ $L^-$-forks over $I'_{\alpha^*}{I_*^{>\alpha_0}}$. Now, as $tp_{L^-}(a_i/a_{<i}I'_{\alpha^*})$ is finitely satisfiable in
$I'_{\alpha^*}$, we conclude that $\phi^-(x,b)$ $L^-$-forks over $a_iI'_{\alpha^*}$ and in particular over $a_i$. Contradiction to $\phi^-(a_i,b)$.\qed\\
 \\Recall the following fact:

\begin{fact}\label{tau extensions}\em [S1, Lemma 2.6]\em
Let $\UU$ be a forking-open set over a set $A$ and let $B\supseteq A$ be any set. Then $\UU$ is
forking-open over $B$.
\end{fact}

\noindent By Claim \ref{claim1} and Fact \ref{tau extensions}, we conclude that $U^T_{\phi^-}$ is a forking-open set over $I_{\FF}$ for $T^-$.\qed

$\\$
\section{An example of a non $L^-$-invariant basic $NFI$-open set}
In [S0], it was observed that for stable theories an invariant set is a basic $NFI$-open set iff it is both $L$-definable over $\emptyset$ and $L^-$-definable with parameters.
Here we show an example of a simple theory with a basic $NFI$-open set that is not even $L^-$-invariant over parameters.\\

Let $L=\{P,R\}$. Let $T_0$ be the $L$-theory that says that $R$ is a symmetric irreflexive binary relation and that $R$ is a complete graph on the unary predicate $P$.
\begin{claim}\label{claim2}
$T_0$ has a model companion $T^*$ and it is axiomatized by the following universal-existential sentences in addition to that of $T_0$:\\
1) For every disjoint finite sets $A,B$ there exists an element $c$ such that $\neg P(c)$ and $R(c,a)$ for all $a\in A$ and such that $\neg R(c,b)$ for all $b\in B$.\\
2) For every disjoint finite sets $A,B$ such that $B\cap P^\CC=\emptyset$ there exists an element $c$ such that $P(c)$ and $R(c,a)$ for all  $a\in A$ and such that $\neg R(c,b)$ for all  $b\in B$.
\end{claim}

\proof We claim that a model $M$ of $T_0$ is existentially closed  for $T_0$ iff $M$ satisfies 1) and 2) of the claim and thus the theory $T^*$ that is obtained from $T_0$ by adding the sentences in 1) and 2) is the model companion of $T_0$. Left to right is immediate. To show the other direction, assume $M$ is a model of $T_0$ that satisfies 1) and 2). Let $\bar a\subseteq M$ be a finite tuple and let $N$ be a model of $T_0$ such that $M$ is a substructure of $N$. Let $\bar b=(b_0,...,b_n)$ be any finite tuple from $N$. It will be sufficient to realize $tp_{qf}^N(\bar b/\bar a)$ (=the quantifier free type of $\bar b$ over $\bar a$ in $N$) in $M$. Indeed, let us construct by induction a sequence $(b'_i \vert i\leq n)\subseteq M$ in the following way. Assume we have constructed $(b'_i \vert i<i_0)$ such that  $tp^M_{qf}(b'_{<i_0}\bar a)=tp^N_{qf}(b_{<i_0}\bar a)$ for some $i_0<n$. We can choose $b'_{i_0}\in M$ such that $tp^M_{qf}(b'_{\leq i_0}\bar a)=tp^N_{qf}(b_{\leq i_0}\bar a)$ in the following way: choose $b'_{i_0}\in M$
such that for all $c\in b_{<i_0}\bar a$ and corresponding $c'\in b'_{<i_0}\bar a$ (i.e. $c,c'$ appears in the same location in the above sequences), we have $R(b_{i_0},c)$ iff $R(b'_{i_0},c')$ and such that $P(b_{<i_0})$ iff $P(b'_{<i_0})$. This can be done using the sentences in 1) in case $\neg P(b_{i_0})$ holds, and using the sentences in 2) in case $P(b_{i_0})$ holds. We can also guarantee that
$b_{i_0}\in b_{<i_0}\bar a$ iff $b'_{i_0}\in b'_{<i_0}\bar a$ (if $b_{i_0}\not\in b_{<i_0}\bar a$, we can guarantee $b'_{i_0}\not\in b'_{<i_0}\bar a$ by choosing  sufficiently many other elements, say $\bar d$, realizing $(\neg P)^M$ and sufficiently many elements in $M$ whose $\{P,R\}$-type over $b'_{<i_0}\bar a$ is $tp_{\{P,R\}}(b'_{i_0}/b'_{<i_0}\bar a)$ but have distinct $R$-types over $\bar d$).\qed

\begin{claim}\label{claim3}
$T^*$ is a complete $\aleph_0$-categorical theory with elimination of quantifiers and $P(x)$ is a complete type over $\emptyset$.
\end{claim}
\proof $\aleph_0$-categoricity is an easy back and forth argument and elimination of quantifiers follows by a similar argument (any isomorphism between finite substructures of the countable model of $T^*$ can be extended to an automophism) and in particular $P(x)$ is a complete type over $\emptyset$.\qed\\

From now on, we work in a highly saturated and highly strongly homogeneous model $\CC$ of $T^*$.

\begin{claim}\label{claim4}
$T^*$ supersimple of $SU$-rank 1. For every $a\in\CC^1$, $P(x)\wedge\neg R(x,a)$ $L$-doesn't fork over $\emptyset$ iff $\neg P(a)$.
\end{claim}

\proof  To prove that $SU(x=x)=1$, note that if $\phi(x,\bar a)$ (where $x$ ia single variable, $\bar a=(a_0,...a_k)$) is not algebraic, then either:\\
i) $P(x)\wedge\bigwedge_{i\leq k} R^{t_i}(x,a_i)\vdash \phi(x,\bar a)$, for $t_i< 2$, where $P(a_i)$ implies $t_i=1$ or\\
ii) $\neg P(x)\wedge\bigwedge_{i\leq k} R^{t_i}(x,a_i)\vdash \phi(x,\bar a)$ for some $t_i< 2$.\\
In each case, using the sentences of 1) and 2), it is easy to conclude that $\phi(x,\bar a)$ doesn't divide over $\emptyset$. The second part is immediate using the sentences in 2) and the definition of $T_0$.\qed

\begin{claim}\label{claim5}
Let $L^-=\{R\}$. Then $P(x)$ is not $L^-$-invariant over parameters and thus $P(x)$ is not an $L^--IF$ type.
\end{claim}

\proof Assume by contradiction that $P(x)$ is $L^-$-invariant over parameters. Then $\bigwedge_{i\leq k} R^{t_i}(x,a_i)\vdash P(x)$ for some $t_i<2$ and some $a_i\in\CC^1$ , but this contradicts the sentences in 1).\qed

\section{$L^-$-invariance of the $NFI$-topology in lovely pairs}
In the other direction we show that in every lovely pair expansion of a theory with the \em wnfcp \em, every complete type over $\emptyset$ is $L^--IF$.\\

\noindent Recall first the basic notions of lovely pairs. Given
$\kappa\geq |T|^{+}$, an elementary pair  $(N,M)$ of models $M\subseteq N$ of a simple theory $T$ is said to
be \em $\kappa$-lovely \em if (i) it has \em the extension property: \em for any $A\subseteq N$ of
cardinality $<\kappa$ and finitary $p(x)\in S(A)$, some nonforking extension of $p(x)$ over $A\cup
M$ is realized in $N$, and (ii) it has \em the coheir property: \em if $p$ as in (i) does not fork
over $M$ then $p(x)$ is realized in $M$. By a \em lovely pair \em (of models of $T$) we mean a
$|T|^{+}$-lovely pair. Let $L_{P}$ be $L$ together with a new unary predicate $P$. Any elementary pair $(N,M)$ of models of $T$ ($M\subseteq N$) can
be considered as an $L_{P}$-structure by taking $M$ to be the interpretation of $P$. A basic property
from [BPV] says that any two lovely pairs of models of $T$ are elementarily equivalent, as
$L_{P}$-structures. So $T_{P}$, the common $L_{P}$-theory of lovely pairs, is complete.
$T$ has the wnfcp if every $\vert T\vert^+$-saturated model of $T_P$ is a lovely pair
(equivalently, for every $\kappa\geq\vert T\vert^+$, any $\kappa$-saturated model of $T_P$ is a $\kappa$-lovely pair).
Every theory with the wnfcp is in particular low (low theories is a subclass of simple theories).
By [BPV, Proposition 6.2], if $T$ has the wnfcp then $T_P$ is simple. Thus, this situation is a special case of our general setting in section 2, where $T_P$ is the given theory ($T$ in the general setting) and $T$ is a reduct ($T^-$ in the general setting). So, in this section we assume $T$ has the wnfcp and we work in a $\lambda$-big model $\MM=(\bar M, P({\bar M}))$ of $T_P$ for some large $\lambda$   (so  $P^\MM=P(\bar M)$). $\nonforkempty{}{}$ will denote independence in $\MM$ and $\minusnonforkempty{}{}$ will denote independence in $\bar M=\MM\vert L$. $dcl^{heq}$ denotes definable closure for hyperimaginries of $\bar M$. Recall the following notation: for  $a\in\bar M^{heq}$, let $a^c=Cb^-(a/P(\bar M))$, where $Cb^-$ denotes the canonical base (as a hyperimaginary element) in the sense of $T$.

\begin{proposition}\label{prop2}
Let  $T$ be a theory with the wnfcp and let $T_P$ be the theory of its lovely pairs. Then any complete type over $\emptyset$ is $L^--IF$ over $\emptyset$.
\end{proposition}

\noindent For proving Proposition \ref{prop2}, we will need the following facts. Recall that for $q(x)\in S(M)$, $Cl(q(x))$ denotes the set of formulas $\phi(x,y)$
without parameters that are represented in $q(x)$ (i.e. such that $\phi(x,a)\in q$ for some tuple $a$ from $M$).

\begin{fact}\label{fact3}$[BPV, Corollary 3.11]$
For every $a,a'$ of the same sort, $tp_{L_P}(a)=tp_{L_P}(a')$ iff $Cl(tp_L(a/P(\bar M)))=Cl(tp_L(a'/P(\bar M)))$.
\end{fact}

\begin{fact} {\em[BPV, Proposition 7.3]\em} \label{fact2}\\
Let $B\subseteq\MM$ and $a$ a tuple from $\MM$. Then\\ $\nonforkempty{a}{B}$ iff $[\minusnonfork{a}{B\cup P(\bar M)}{P(\bar M)}$ and $\minusnonforkempty{a^c}{B^c}]$.
\end{fact}

\begin{fact}\label{fact1}$[S0, Proposition 3.1]$
For $T_P$ and the reduct $T$, we have: for every finite tuple of variables $x$, $\tilde\Gamma_{x}=(x=x)$, namely the greatest universal transducer in the variables $x$ is  $(x=x)$.
\end{fact}

\begin{claim}\label{claim6}
Let  $p(x)\in S_x(T_P)$. Let $a\models p$ and let $q=tp_L(a,a^c)$. If $\hat a\in \MM$ is such that for some $a'\in bdd(P(\bar M))$, $(\hat a,a')\models q$ and $\minusnonfork{\hat a}{P(\bar M)}{a'}$.
Then $tp_{L_P}(\hat a)=p$. Clearly, $a'$ is a canonical base of $tp_L(\hat a/P(\bar M))$.
\end{claim}

\proof By the assumption, for every $\psi^-( x,y)\in L$, $\psi^-(\hat a,y)$ $L$-doesn't fork over $a'$ iff $\psi^-(\hat a,y)$ $L$-doesn't fork over $P(\bar M)$ iff $\psi^-(\hat a,y)$ is realized in
$P(\bar M)$. Likewise, as $\minusnonfork{a}{P(\bar M)}{a^c}$, we conclude that $\psi^-(a,y)$ $L$-doesn't fork over $a^c$ iff $\psi^-(a,y)$ is realized in $P(\bar M)$.
Thus $Cl(tp_L(a/P(\bar M)))=Cl(tp_L(\hat a/P(\bar M)))$ and so by Fact \ref{fact3},  $tp_{L_P}(\hat a)=tp_{L_P}(a)=p$. \qed\\

\noindent \textbf{Proof of Proposition \ref{prop2}} Let $p(x)\in S_x(T_P)$ and $\phi^-(x,y)\in L$. For $a\models p$, let $q=tp_L(a,a^c)$ (clearly $a^c\in dcl^{heq}(P(\bar M)$). Then, clearly $q$ depends only on $p$. It will be sufficient to prove the following ($p^-(x)$ denotes $p(x)\vert L$).
\begin{claim}\label{claim7}
For all $b$, we have:
$$b\in U_{p,\phi^-}\mbox{ iff }\ p^-(x)\wedge\phi^-(x,b)\ \mbox{L-doesn't\ fork\ over}\ \emptyset.$$
\end{claim}


\proof Assume $b\in U_{p,\phi^-}$. Then there exists $a\models p$ such that $\nonforkempty{a}{b}$ and $\phi^-(a,b)$.  By Fact \ref{fact1}, $\minusnonforkempty{a}{b}$. To prove the other direction, assume $p^-(x)\wedge\phi^-(x,b)\ L\mbox{-doesn't\ fork\ over}\ \emptyset$. Then,  there exists $a\models p^-(x)\wedge\phi^-(x,b)$ such that $\minusnonforkempty{a}{b}$.
Hence, there exists $a'$ such that $\minusnonforkempty{aa'}{b}$ and $(a,a')\models q$. By the extension property, we may assume $\minusnonforkempty{aa'}{bP(\bar M)}$. Then $\minusnonfork{aa'}{bb^c}{P(\bar M)}$ and in particular, $tp_L(aa'/bb^c)$ doesn't fork over $P(\bar M)$. By the coheir property, there are $a_0,a'_0\in dcl^{heq}(P(\bar M))$ such that $a_0a'_0$ realize $tp_L(aa'/bb^c)$. So, clearly $\minusnonforkempty{a_0a'_0}{bb^c}$ and $\phi^-(a_0,b)$.
In particular, $\minusnonfork{a_0}{bb^c}{a'_0}$. By the extension property, there exists $a^*_0\in\MM$ that realizes an $L$-non-forking extension of $tp_L(a_0/a'_0bb^c)$ over $P(\bar M)a'_0bb^c$.
In particular, $\minusnonfork{a^*_0}{P(\bar M)}{a'_0}$. By Claim \ref{claim6}, $tp_{L_P}(a^*_0)=p$. Now, as $\minusnonfork{a^*_0}{b}{P(\bar M)}$, and $a'_0$ is a canonical base of $tp_L(a^*_0/P(\bar M))$ and $\minusnonforkempty{a'_0}{b^c}$ we conclude by Fact \ref{fact2} that $\nonforkempty{a^*_0}{b}$. As clearly $\phi^-(a^*_0,b)$ we conclude that $p(x)\wedge\phi^-(x,b)$ $L_{P}$-doesn't fork over $\emptyset$.\qed

\noindent Ziv Shami, E-mail address: zivsh@ariel.ac.il.\\
Dept. of Mathematics\\
Ariel University\\
Samaria, Ariel 44873\\
Israel.\\

\end{document}